\documentclass[11pt]{amsart}

\usepackage[margin=0.85in]{geometry}
\usepackage{amsmath,amssymb,amsthm,mathtools}
\usepackage{array}
\usepackage{booktabs}

\newtheorem{theorem}{Theorem}[section]
\newtheorem{lemma}[theorem]{Lemma}
\newtheorem{proposition}[theorem]{Proposition}

\theoremstyle{remark}
\newtheorem{remark}[theorem]{Remark}

\newcommand{\F}{\mathbb F}

\newcommand{\Alt}{\operatorname{Alt}}
\newcommand{\per}{\operatorname{per}}

\title{On Distinct Powers in Finite Groups}

\author[D. Y. Ma]{Dingying Ma}
\address{College of Mathematics and Statistics, Yili Normal University, Yining 835000, China}
\email{mady0229@hotmail.com}


\begin{document}

\begin{abstract}
Sun conjectured that if a group $G$ has no elements of orders $2,3,\ldots,n+1$, then every $n$-element subset $A$ can be ordered as $a_1,\ldots,a_n$ so that $a_1,a_2^2,\ldots,a_n^n$ are pairwise distinct. We prove the conjecture for finite groups when $n=4,5,6$. The proofs combine collision-counting arguments with reductions to finite-field assignment problems; the remaining finite cases for $n=6$ are verified exhaustively. We also prove the conjecture for every $n$ when no four elements of $A$ generate the same cyclic subgroup.
\end{abstract}

\maketitle

\section{Introduction}

Sun introduced this conjecture in \cite[Conjecture 4.1]{Sun}. It states that if a group $G$ has no elements of orders $2,3,\ldots,n+1$, then every $n$-element subset $A$ of $G$ can be ordered as $a_1,\ldots,a_n$ so that
\[
a_1,\ a_2^2,\ldots,a_n^n
\]
are pairwise distinct. The conjecture is also Problem 20.100 in the Kourovka Notebook \cite[Problem 20.100]{KM}, where it is noted that the conjecture is known for $n\leq3$. Sun proved these cases and, for every $n$, the case of torsion-free abelian groups \cite[Theorem 1.4]{Sun}. He also reported that the conjecture had been checked by computer for cyclic groups of bounded order when $n=4,\ldots,9$ \cite[Remark 4.1(b)]{Sun}. Thus $n=4$ is the first case not covered by the known general results.

Besides the cases $n=4,5,6$ for finite groups, we prove a result valid for every $n$. For $x,y\in A$, write $x\sim y$ if $\langle x\rangle=\langle y\rangle$. Thus the equivalence classes consist of those elements of $A$ that generate the same cyclic subgroup.

\begin{theorem}\label{thm:class3}
Let $n$ be a positive integer, let $G$ be a finite group with no elements of orders $2,3,\ldots,n+1$, and let $A\subseteq G$ with $|A|=n$. Suppose that every equivalence class of $\sim$ has size at most $3$. Then the elements of $A$ can be ordered as $a_1,\ldots,a_n$ so that
\[
a_1,\ a_2^2,\ldots,a_n^n
\]
are pairwise distinct.
\end{theorem}

\begin{theorem}\label{thm:n4}
Let $G$ be a finite group with no elements of orders $2,3,4,5$. If $A\subseteq G$ and $|A|=4$, then the elements of $A$ can be ordered as $a_1,\ldots,a_4$ so that
\[
a_1,\ a_2^2,\ a_3^3,\ a_4^4
\]
are pairwise distinct.
\end{theorem}

\begin{theorem}\label{thm:main}
Let $G$ be a finite group with no elements of orders $2,3,4,5,6$. If $A\subseteq G$ and $|A|=5$, then the elements of $A$ can be ordered as $a_1,\ldots,a_5$ so that
\[
a_1,\ a_2^2,\ a_3^3,\ a_4^4,\ a_5^5
\]
are pairwise distinct.
\end{theorem}

\begin{theorem}\label{thm:n6}
Let $G$ be a finite group with no elements of orders $2,3,4,5,6,7$. If $A\subseteq G$ and $|A|=6$, then the elements of $A$ can be ordered as $a_1,\ldots,a_6$ so that
\[
a_1,\ a_2^2,\ a_3^3,\ a_4^4,\ a_5^5,\ a_6^6
\]
are pairwise distinct.
\end{theorem}

We begin with Theorem~\ref{thm:class3}, where the elements of $A$ are divided according to the cyclic subgroups they generate. We then treat four elements. For five elements, the only remaining case consists of five generators of one cyclic subgroup, and it reduces to a question over $\F_p$ for primes $p\geq7$. The proof for six elements starts in the same way. Its main new ingredient is a weighted graph whose edges are determined by the ratios $2$, $3$, and their inverses. Except for eighteen primes $p\geq37$, the relevant part of this graph has the same form as the corresponding graph over the rational numbers. Those eighteen primes, together with the seven smaller primes, are treated individually in the appendix.

\section{When no four elements generate the same cyclic subgroup}

We prove Theorem~\ref{thm:class3}. We shall use the following three-term observation.

\begin{lemma}\label{lem:three-block}
Let $p$ be a prime, let $1\leq r$ and $r+2<p$, and suppose that $p\nmid 6(r+1)$. For any $b_1,b_2,b_3\in\F_p^\times$, there is a permutation $\pi\in S_3$ such that
\[
r b_{\pi(1)},\qquad (r+1)b_{\pi(2)},\qquad (r+2)b_{\pi(3)}
\]
are pairwise distinct.
\end{lemma}

\begin{proof}
First suppose that $b_1,b_2,b_3$ are distinct. Put
\[
P(X_1,X_2,X_3)=
\prod_{1\leq i<j\leq3}(a_iX_i-a_jX_j),
\qquad (a_1,a_2,a_3)=(r,r+1,r+2).
\]
For a polynomial $Q(X_1,X_2,X_3)$, define
\[
\Alt(Q)=\sum_{\sigma\in S_3}\operatorname{sgn}(\sigma)
Q(X_{\sigma(1)},X_{\sigma(2)},X_{\sigma(3)}).
\]
A direct expansion gives
\[
\Alt(P)=-6(r+1)^3
\prod_{1\leq i<j\leq3}(X_j-X_i).
\]
The right side is nonzero at $(b_1,b_2,b_3)$. Hence at least one of the six values
$P(b_{\sigma(1)},b_{\sigma(2)},b_{\sigma(3)})$ is nonzero, and the required ordering follows.

Now suppose that exactly two of the $b_i$ are equal. Multiplying all three entries by a common nonzero scalar does not affect whether the three products are pairwise distinct. We may therefore write the multiset as $\{1,1,q\}$ with $q\ne1$. If every placement of $q$ failed, then for every
\[
a\in B:=\{r,r+1,r+2\}
\]
we would have $qa\in B\setminus\{a\}$. Thus multiplication by $q$ would permute $B$ without a fixed point. It would therefore be a $3$-cycle, and comparison of the sums of the elements of $B$ and $qB$ would give
\[
3(r+1)=3q(r+1)
\]
in $\F_p$. Since $p\nmid3(r+1)$, this forces $q=1$, a contradiction. If all three entries are equal, the assertion follows immediately because $r,r+1,r+2$ are distinct modulo $p$.
\end{proof}

\begin{proof}[Proof of Theorem~\ref{thm:class3}]
Every prime divisor of the order of a nonidentity element of $G$ is greater than $n+1$ by the hypothesis. If $x^r=y^s$ with $x,y\in A$ and $1\le r,s\le n$, then $r$ and $s$ are coprime to $|x|$ and $|y|$, respectively, and hence $\langle x\rangle=\langle x^r\rangle=\langle y^s\rangle=\langle y\rangle$. Thus powers coming from different $\sim$-classes cannot coincide. Suppose that there are $s$ classes of size $2$, $t$ classes of size $3$, and $u$ singleton classes. Thus
\[
2s+3t+u=n.
\]
Give the two-element classes, in any order, the exponent pairs
\[
\{1,n\},\ \{2,n-1\},\ldots,\{s,n+1-s\}.
\]
Consider one such class $\{x,y\}$ with exponent pair $\{r,n+1-r\}$. If both assignments failed, then
\[
x^r=y^{n+1-r},\qquad y^r=x^{n+1-r}.
\]
Hence the common order of $x$ and $y$ would divide
\[
(n+1-r)^2-r^2=(n+1)(n+1-2r).
\]
Here $r\le s\le n/2$, so $n+1-2r>0$. Every prime divisor of this positive integer is at most $n+1$, whereas every prime divisor of the common order is greater than $n+1$. This is impossible. Thus each two-element class can be arranged independently.

The unused exponents form the interval
\[
s+1,s+2,\ldots,n-s,
\]
which has length $n-2s=3t+u$. Give the $t$ three-element classes consecutive triples from the beginning of this interval, and give the remaining $u$ exponents to the singleton classes. Consider a three-element class and one of these triples
\[
\{r,r+1,r+2\}.
\]
Let its elements generate the cyclic group $\langle g\rangle$ of order $m$, and write them as $g^{v_1},g^{v_2},g^{v_3}$, where each $v_i$ is coprime to $m$. Choose a prime $p\mid m$. Then $p>n+1$, so in particular $p>r+2$ and $p\nmid6(r+1)$. Reducing the $v_i$ modulo $p$ and applying Lemma~\ref{lem:three-block}, we may assign the three exponents so that the products of the assigned exponents with the corresponding $v_i$ are pairwise distinct modulo $p$. They are therefore pairwise distinct modulo $m$, and the resulting three powers of $g$ are pairwise distinct.

We have arranged each $\sim$-class so that its selected powers are distinct, and powers from different classes cannot coincide by the observation at the beginning of the proof. This proves the theorem.
\end{proof}

\begin{remark}\label{rem:four-block}
The consecutive-block argument does not extend directly from three elements to four. For example, in $\F_{11}$ take the coefficient block $\{4,5,6,7\}$ and the multiset $\{1,1,1,-1\}$. Wherever $-1$ is placed, two of the four products coincide. Thus a four-element generator class requires a different choice of exponents or a different argument. This example only shows a limitation of the consecutive-block proof; it does not show that the hypothesis in Theorem~\ref{thm:class3} is sharp. It remains open whether an analogous general theorem holds when the equivalence classes have size at most $4$.
\end{remark}

\section{Counting possible coincidences}

The same elementary count will be used for four, five, and six elements, so we state it once here. Let $B\subseteq\mathbb Z$ be a finite set with $|B|=t$, let $p$ be a prime, and let $q\in\F_p^\times$. Define
\[
w_{B,p}(q)=\#\{(r,s)\in B^2:r\ne s,\ r\equiv qs\pmod p\}.
\]
Thus $w_{B,p}(q)$ counts the ordered pairs of distinct coefficients in $B$ whose ratio is $q$ modulo $p$. For $c_1,\ldots,c_t\in\F_p^\times$, put
\[
W_B(c_1,\ldots,c_t)=\sum_{1\le i<j\le t}w_{B,p}(c_j/c_i).
\]

\begin{lemma}\label{lem:collision-count}
With the notation above, if
\[
W_B(c_1,\ldots,c_t)<t(t-1),
\]
then the entries $c_1,\ldots,c_t$ can be assigned bijectively to the coefficients in $B$ so that the resulting products are pairwise distinct.
\end{lemma}

\begin{proof}
There are $t!$ assignments. For a fixed unordered pair of labeled entries $\{c_i,c_j\}$, the two resulting products are equal precisely when their assigned coefficients $(r,s)$ satisfy $r c_i=s c_j$. By the definition of $w_{B,p}$, there are $w_{B,p}(c_j/c_i)$ ordered coefficient pairs that produce such an equality, and each occurs in exactly $(t-2)!$ assignments. Hence, if we sum over all assignments and all unordered pairs of entries, the total number of equalities is
\[
(t-2)!\,W_B(c_1,\ldots,c_t).
\]
Therefore the average number of equal pairs in an assignment is
\begin{equation}\label{eq:expectation-general}
\frac{W_B(c_1,\ldots,c_t)}{t(t-1)}.
\end{equation}
If $W_B<t(t-1)$, this average is less than $1$, so at least one assignment gives pairwise distinct products.
\end{proof}

\section{Four elements}

Set $B_4=\{1,2,3,4\}$ and write $w_{4,p}(q)=w_{B_4,p}(q)$.

\begin{lemma}\label{lem:field4}
Let $b_1,b_2,b_3,b_4\in\F_p^\times$.
\begin{enumerate}
\item If $p\ge11$, there is a permutation $\pi\in S_4$ such that
\[
b_{\pi(1)},\quad 2b_{\pi(2)},\quad 3b_{\pi(3)},\quad4b_{\pi(4)}
\]
are pairwise distinct.
\item If $p=7$, the same conclusion holds unless, up to multiplication of all four entries by a nonzero scalar, the multiset is
\[
\{1,1,2,2\}.
\]
\end{enumerate}
\end{lemma}

\begin{proof}
If every assignment were bad, Lemma~\ref{lem:collision-count} would give
\begin{equation}\label{eq:weight4}
W_{B_4}(b_1,b_2,b_3,b_4)=\sum_{1\le i<j\le4}w_{4,p}(b_j/b_i)\ge12.
\end{equation}
For $p\ge11$, one has $w_{4,p}(q)\le2$. Thus equality must hold in \eqref{eq:weight4}, the four entries must be distinct, and every ratio between two entries must have weight $2$. For $p\ge17$, the only such ratios are $2$ and $1/2$, which is already impossible for three distinct entries. For $p=11$ and $13$, the weight-$2$ ratio sets are, respectively,
\[
\{2,3,4,6\},\qquad \{2,4,7,10\}.
\]
After scaling one entry to $1$, the displayed sets show that no three further elements can have all their pairwise ratios in the same set. Hence no four-element counterexample exists.

For $p=7$, the weights are
\[
w_{4,7}(2)=w_{4,7}(4)=3,
\qquad
w_{4,7}(3)=w_{4,7}(5)=w_{4,7}(6)=2.
\]
The multiplicity types $4$ and $3+1$ have $W_{B_4}\le9$, so Lemma~\ref{lem:collision-count} applies. In type $2+2$, after scaling the multiset is $\{1,1,q,q\}$ and
\[
W_{B_4}=4w_{4,7}(q).
\]
Thus this type is good for $q=3,5,6$. For $q=2$ or $4$, each of the six possible placements of the two equal entries gives an equality, giving the exceptional scalar orbit $\{1,1,2,2\}$.

It remains to justify the types $2+1+1$ and $1+1+1+1$. In the first type, scale the repeated entry to $1$. For the three choices $\{1,1,3,5\}$, $\{1,1,3,6\}$, and $\{1,1,5,6\}$ one has $W_{B_4}=11$, so Lemma~\ref{lem:collision-count} applies. For the other seven possibilities the following table gives, in the second column, an ordering assigned to the coefficients $1,2,3,4$; the resulting products are immediately seen to be distinct modulo $7$.
\[
\begin{array}{c|c}
\text{multiset}&\text{ordering}\\ \hline
\{1,1,2,3\}&(1,1,2,3)\\
\{1,1,2,4\}&(4,1,1,2)\\
\{1,1,2,5\}&(1,5,2,1)\\
\{1,1,2,6\}&(2,6,1,1)\\
\{1,1,3,4\}&(1,3,1,4)\\
\{1,1,4,5\}&(1,5,4,1)\\
\{1,1,4,6\}&(6,4,1,1)
\end{array}
\]
Finally, up to scaling, the four-distinct multisets are
\[
\{1,2,3,4\},\qquad \{1,2,3,5\},\qquad \{1,2,5,6\},
\]
and the orderings
\[
(2,4,1,3),\qquad (2,3,5,1),\qquad (2,6,5,1)
\]
respectively give four distinct products modulo $7$. Hence the only bad multisets are the scalar multiples of $\{1,1,2,2\}$. Equivalently, the six bad multisets are
\[
\{1,1,2,2\},\ \{1,1,4,4\},\ \{2,2,4,4\},\
\{3,3,5,5\},\ \{3,3,6,6\},\ \{5,5,6,6\}.
\]
\end{proof}

\begin{proof}[Proof of Theorem~\ref{thm:n4}]
If $x^r=y^s$ with $1\le r,s\le4$, then $\langle x\rangle=\langle y\rangle$, since every prime divisor of an element order is at least $7$. Partition $A$ by equality of the generated cyclic subgroup. Different classes can never produce equal selected powers.

If the class sizes are $3+1$, use Sun's $n=3$ case \cite[Theorem 1.4(ii)]{Sun} on the three-element class with exponents $1,2,3$, and give exponent $4$ to the singleton. If the class sizes are $2+2$, use exponent pairs $\{1,3\}$ and $\{2,4\}$. If both assignments failed on a pair with exponents $r<s$, then the common element order would divide $s^2-r^2$; here the two values are $8$ and $12$, impossible. If the class sizes are $2+1+1$, use exponents $\{1,2\}$ on the two-element class and give exponents $3,4$ to the two singletons. If both assignments on the pair failed, the same calculation as above would force the common element order to divide $2^2-1^2=3$, which is impossible. If all four classes are singletons, no two selected powers can coincide. Thus it remains only to consider four elements generating the same cyclic subgroup.

Write
\[
A=\{g^{u_1},g^{u_2},g^{u_3},g^{u_4}\},\qquad |g|=m,
\]
where each $u_i$ is a unit modulo $m$. If $m$ has a prime divisor $p\ge11$, reduce the $u_i$ modulo $p$ and apply Lemma~\ref{lem:field4}(1). Distinctness modulo $p$ implies distinctness modulo $m$.

It remains to consider $m=7^a$. Reduce the $u_i$ modulo $7$. If their multiset is not exceptional in Lemma~\ref{lem:field4}(2), we are done. Otherwise, after multiplying all four $u_i$ by one common unit, we may write their residues modulo $7$ as
\[
u_1,u_2\equiv1,\qquad v_1,v_2\equiv2.
\]
This normalization is legitimate: replacing $g$ by another generator of $\langle g\rangle$ multiplies all exponent residues modulo $7$ by the same nonzero scalar.
There are indices $i,j$ with
\[
v_j\ne2u_i\pmod{7^a}.
\]
Indeed, if $v_j=2u_i$ for every $i,j$, then $v_1=v_2$ modulo $7^a$, contrary to the distinctness of the four elements.
Assign the four entries in the order
\[
(v_j,u_i,u_{3-i},v_{3-j})
\]
to the exponents $1,2,3,4$. Modulo $7$ the four products are
\[
2,\quad2,\quad3,\quad1,
\]
so the only possible collision modulo $7^a$ is between the first two. Our choice of $i,j$ excludes that collision. Hence the four selected powers are pairwise distinct.
\end{proof}

\section{Reducing to a single cyclic subgroup}

We first note when two selected powers can be equal.

\begin{lemma}\label{lem:same-cyclic}
Let $G$ be finite and suppose that $G$ has no elements of orders $2,3,4,5,6$. If $x,y\in G$ and
\[
x^r=y^s
\]
for some $1\leq r,s\leq 5$, then
\[
\langle x\rangle=\langle y\rangle.
\]
\end{lemma}

\begin{proof}
Every prime divisor of $|x|$ or $|y|$ is at least $7$. Hence $r$ and $s$ are coprime to $|x|$ and $|y|$, respectively. Thus
\[
\langle x\rangle=\langle x^r\rangle
=\langle y^s\rangle=\langle y\rangle.
\]
\end{proof}

We shall also use the following two-point observation.

\begin{lemma}\label{lem:pair}
Let $x\neq y$ generate the same finite cyclic group, and let $1\leq r<s$. If both assignments of the exponents $r,s$ to $x,y$ give a collision, then
\[
x^{s^2-r^2}=1.
\]
\end{lemma}

\begin{proof}
The two failures are
\[
x^r=y^s,\qquad y^r=x^s.
\]
Raising the first equality to the $r$th power and using the second gives
\[
x^{r^2}=y^{rs}=x^{s^2}.
\]
\end{proof}

\begin{proposition}\label{prop:cyclic-reduction}
To prove Theorem~\ref{thm:main}, it is enough to consider the case in which all five elements of $A$ generate the same cyclic subgroup.
\end{proposition}

\begin{proof}
Partition $A$ by the relation
\[
x\sim y \quad\Longleftrightarrow\quad \langle x\rangle=\langle y\rangle.
\]
By Lemma~\ref{lem:same-cyclic}, powers attached to elements from two different classes can never be equal. Hence the exponent sets may be chosen separately on the different classes.

The possible nontrivial partitions of $5$ are dealt with as follows:
\[
\begin{array}{c|c}
\text{class sizes}&\text{assigned exponent sets}\\
\hline
4+1&\{1,2,3,4\}+\{5\},\\
3+2&\{1,2,3\}+\{4,5\},\\
3+1+1&\{1,2,3\}+\{4\}+\{5\},\\
2+2+1&\{1,3\}+\{2,4\}+\{5\},\\
2+1+1+1&\{1,2\}+\{3\}+\{4\}+\{5\}.
\end{array}
\]
For classes of size $3$ we use the known case $n=3$, and for a class of size $4$ we use Theorem~\ref{thm:n4}. For a two-element class, Lemma~\ref{lem:pair} applies. The relevant values of $s^2-r^2$ are
\[
9,\qquad 8,\qquad 12,\qquad 3,
\]
for the exponent pairs $\{4,5\}$, $\{1,3\}$, $\{2,4\}$, and $\{1,2\}$, respectively. Each has all prime divisors in $\{2,3\}$, whereas every nonidentity element of $G$ has order with all prime divisors at least $7$. Hence at least one assignment works on every two-element class.
\end{proof}

Suppose now that
\[
A=\{g^{u_1},\ldots,g^{u_5}\}\subseteq \langle g\rangle,
\qquad |g|=m,
\]
and every member of $A$ generates $\langle g\rangle$. Then $(u_i,m)=1$ for every $i$. Choose a prime $p\mid m$. Necessarily $p\geq 7$, and each $u_i$ is nonzero modulo $p$. If two integer multiples of the $u_i$ are distinct modulo $p$, then they are distinct modulo $m$, since congruence modulo $m$ would imply congruence modulo $p$. The following lemma will therefore finish the proof.

\section{Five elements over a finite field}

\begin{lemma}\label{lem:field}
Let $p\geq 7$ be prime and let $b_1,\ldots,b_5\in\F_p^\times$, with repetitions allowed. Then there is a permutation $\pi\in S_5$ such that
\[
b_{\pi(1)},\quad 2b_{\pi(2)},\quad 3b_{\pi(3)},\quad
4b_{\pi(4)},\quad 5b_{\pi(5)}
\]
are pairwise distinct.
\end{lemma}

For a multiset, its \emph{support} is the set of its distinct entries. We treat $p=7$, $p=11$, and $p\geq13$ separately.

\subsection{The prime $7$}

Fix a generator of $\F_7^\times$ and take discrete logarithms. Since the coefficients $1,2,3,4,5$ represent all nonzero residues except $-1$, their logarithms form
\[
R=C_6\setminus\{3\}=\{0,1,2,4,5\}.
\]
Thus it is enough to prove the following additive statement: for every five-term multiset $S\subset C_6$, its elements can be assigned to the rows $R$ so that the five sums are distinct.

Translation and negation of $S$ preserve this property. If $S$ has support $1$, the assertion is immediate. If
\[
S=\{0,0,0,0,d\},\qquad d\neq0,
\]
place $d$ in row $3-d$ and the zeros elsewhere, where the row index is taken in $C_6$. The resulting sums are $3$ together with the other four rows.

If
\[
S=\{0,0,0,d,d\},
\]
and $d\neq3$, place the two copies of $d$ in rows
\[
3-d,\qquad 3-2d.
\]
The resulting sums are $3$ and $3-d$, and the other three rows are unchanged. For $d=3$, place the two copies of $3$ in rows $1$ and $4$.

Now suppose
\[
S=\{0,0,0,a,b\},\qquad a,b\neq0,\quad a\neq b.
\]
If $a+b\neq0$, put $a$ in row $x=3-a$ and $b$ in row $y=3-a-b$. Then
\[
x+a=3,\qquad y+b=x,
\]
and the other three rows are unchanged. If $a+b=0$, then, up to negation, only the following possibility remains:
\[
\{0,0,0,1,5\},\qquad \{0,0,0,2,4\}.
\]
On the rows $(0,1,2,4,5)$, the assignments
\[
(0,0,0,1,5),\qquad (2,0,4,0,0)
\]
give the distinct sums
\[
(0,1,2,5,4),\qquad (2,1,0,4,5),
\]
respectively.

It remains to consider the multiplicity types $2+2+1$, $2+1+1+1$, and $1+1+1+1+1$. Translate one entry to $0$, sort the multiset, and identify a multiset with its negative. These two symmetries reduce the possibilities to the following thirteen representatives. In the first column, a string such as $00112$ denotes the multiset $\{0,0,1,1,2\}$. The table gives an assignment to the rows $(0,1,2,4,5)$ and the resulting sums in $C_6$; the table therefore settles every remaining case.

\[
\begin{array}{c|c|c}
S&\text{assignment}&\text{resulting sums}\\
\hline
00112&(0,1,1,0,2)&(0,2,3,4,1)\\
00113&(3,0,0,1,1)&(3,1,2,5,0)\\
00122&(2,0,1,2,0)&(2,1,3,0,5)\\
00123&(0,3,0,1,2)&(0,4,2,5,1)\\
00124&(4,0,1,2,0)&(4,1,3,0,5)\\
00125&(2,5,1,0,0)&(2,0,3,4,5)\\
00133&(3,0,3,0,1)&(3,1,5,4,0)\\
00134&(3,4,0,0,1)&(3,5,2,4,0)\\
00135&(1,3,0,5,0)&(1,4,2,3,5)\\
00144&(4,0,1,4,0)&(4,1,3,2,5)\\
00224&(2,2,4,0,0)&(2,3,0,4,5)\\
00234&(2,0,3,0,4)&(2,1,5,4,3)\\
01234&(0,3,1,4,2)&(0,4,3,2,1)
\end{array}
\]
This proves Lemma~\ref{lem:field} for $p=7$.

\subsection{The prime $11$}

Since $4$ has order $5$ modulo $11$, we have the direct product decomposition
\[
\F_{11}^\times=\{\pm1\}\times\langle4\rangle.
\]
Thus every element of $\F_{11}^\times$ has a unique expression
\[
\varepsilon 4^q,
\qquad \varepsilon\in\{\pm1\},\quad q\in C_5.
\]
The five coefficients $1,2,3,4,5$ represent the five distinct powers of $4$ modulo sign; explicitly,
\[
1=4^0,\qquad 2=-4^3,\qquad 3=4^4,\qquad 4=4,\qquad 5=4^2.
\]
After translating the row labels, we may therefore take the rows to be $0,1,2,3,4$, with row $0$ corresponding to the unique negative coefficient. If $(\varepsilon,q)$ is placed in row $r$, its quotient class is $r+q$, up to a fixed translation, and its sign is reversed exactly when $r=0$. A fixed translation of all resulting quotient classes does not affect whether they are distinct.

The affine changes
\[
q\mapsto aq+b,\qquad r\mapsto ar,
\qquad a\in C_5^\times,
\]
preserve distinctness of the resulting quotient classes and keep the distinguished row $0$ fixed. Translating one quotient entry to $0$, sorting, and then applying the four possible nonzero multipliers $a$ gives nine orbits of five-term quotient multisets. Representatives for all nine orbits appear in the two tables below.

The four zero-sum types have complete transversals. The following table lists an assignment to rows $0,1,2,3,4$ and the five resulting quotient classes.
\[
\begin{array}{c|c|c}
Q&\text{assignment}&\text{resulting quotient classes}\\
\hline
00000&(0,0,0,0,0)&(0,1,2,3,4)\\
00014&(1,4,0,0,0)&(1,0,2,3,4)\\
00113&(3,0,0,1,1)&(3,1,2,4,0)\\
01234&(1,2,3,4,0)&(1,3,0,2,4)
\end{array}
\]
Thus the signs do not matter in these cases.

The five nonzero-sum types are
\[
00001,\qquad00011,\qquad00012,\qquad00112,\qquad00123.
\]
For each type we use two assignments. In the first, one member of a specified doubled pair lies in the negative row $0$; in the second, neither member does. The assignments may be chosen so that the same two labeled entries form the doubled quotient class.
\[
\begin{array}{c|c|c|c}
Q&T_1&T_2&\text{doubled quotient values}\\
\hline
00001&00001&00010&\{0,1\}\\
00011&00011&00110&\{0,1\}\\
00012&00120&00012&\{0,2\}\\
00112&01120&00112&\{0,2\}\\
00123&00123&01230&\{0,2\}
\end{array}
\]
Here the strings give the entries placed in rows $0,1,2,3,4$. A direct addition modulo $5$ shows that each assignment has quotient multiplicities $2,1,1,1,0$ and has the stated doubled pair.

Choose labeled copies realizing the doubled pair. If their original signs are equal, use $T_1$. Exactly one of the two lies in row $0$, so their final signs are opposite. If their original signs are opposite, use $T_2$. Neither lies in row $0$, so their final signs remain opposite. In either case the two elements in the repeated quotient class are distinct, while all other quotient classes occur only once. This proves Lemma~\ref{lem:field} for $p=11$.

\subsection{Primes $p\geq13$}

Set $B_5=\{1,2,3,4,5\}$ and write $w_{5,p}(q)=w_{B_5,p}(q)$. If every assignment of the five entries to the coefficients in $B_5$ were bad, Lemma~\ref{lem:collision-count} would give
\begin{equation}\label{eq:weight}
W_{B_5}(b_1,\ldots,b_5)=\sum_{1\leq i<j\leq5}w_{5,p}(b_j/b_i)\geq20.
\end{equation}
For all primes $p\geq13$ one has $w_{5,p}(q)\leq2$. In particular, if two of the $b_i$ are equal, then one pair contributes $w_{5,p}(1)=0$, so the left side of \eqref{eq:weight} is at most $18$. Thus a counterexample must have five distinct entries, and equality in \eqref{eq:weight} forces every pair to have weight $2$.

For $p\geq29$, two fractions $r/s$ and $r'/s'$ with $1\leq r,s,r',s'\leq5$ can be congruent modulo $p$ only when they are equal as rational numbers, since
\[
|rs'-r's|\leq24<p.
\]
The only repeated ratios are
\[
2=\frac21=\frac42,
\qquad
\frac12=\frac12=\frac24.
\]
Hence equality in \eqref{eq:weight} would force every ratio between two of the five entries to be $2$ or $1/2$, which is impossible already for three distinct entries.

It remains to consider $p=13,17,19,23$. For $p=17,19,23$, put
\[
D_p=\{q\in\F_p^\times:w_{5,p}(q)=2\}.
\]
The relevant sets are
\[
\begin{array}{c|c|c}
p&D_p&\text{largest clique containing }1\\
\hline
17&\{2,4,5,7,9,11,13,14\}&5\\
19&\{2,4,5,8,10,12\}&4\\
23&\{2,5,12,14\}&3
\end{array}
\]
Here vertices $x,y\in\F_p^\times$ are adjacent when $x/y\in D_p$. The displayed clique bounds follow by normalizing one vertex to $1$ and successively intersecting the neighborhoods $xD_p$ of the remaining vertices. Thus $p=19,23$ cannot satisfy \eqref{eq:weight}. For $p=17$, the same neighborhood calculation leaves a single normalized $5$-clique, so every $5$-clique is a scalar multiple of
\[
\{1,2,4,5,11\}.
\]
For this set, assigning
\[
(1,2,4,11,5)
\]
to the coefficients $1,2,3,4,5$ gives
\[
1,4,12,10,8,
\]
which are distinct modulo $17$.

Finally let $p=13$. As above, repeated entries cannot give a counterexample, so assume that $b_1,\ldots,b_5$ are distinct. Set
\[
F(X_1,\ldots,X_5)=\prod_{1\leq r<s\leq5}(rX_r-sX_s).
\]
For a polynomial $P(X_1,\ldots,X_5)$, write
\[
\Alt(P)=\sum_{\sigma\in S_5}\operatorname{sgn}(\sigma)
P(X_{\sigma(1)},\ldots,X_{\sigma(5)}).
\]
For a square matrix $M=(m_{ij})$, its permanent is
\[
\per(M)=\sum_{\sigma\in S_5}\prod_{i=1}^5 m_{i,\sigma(i)}.
\]
The permanent below is the same one that occurs in Sun's coefficient identity \cite[Lemma 2.2]{Sun}, which he uses in the torsion-free abelian case. Since $\Alt(F)$ is an alternating polynomial of degree $10$, we have
\[
\Alt(F)=C_5\prod_{i<j}(X_j-X_i)
\]
for an integer $C_5$. Comparing the coefficient of $X_1^0X_2^1X_3^2X_4^3X_5^4$ gives
\[
C_5=\per(i^{j-1})_{1\leq i,j\leq5}=8527540\equiv8\pmod{13}.
\]
Hence
\[
\Alt(F)(b_1,\ldots,b_5)\neq0.
\]
At least one summand in the antisymmetrization is nonzero, and the corresponding permutation gives five distinct products. This proves Lemma~\ref{lem:field} for every $p\geq13$.

\section{Five elements}

\begin{proof}[Proof of Theorem~\ref{thm:main}]
By Proposition~\ref{prop:cyclic-reduction}, we may assume
\[
A=\{g^{u_1},\ldots,g^{u_5}\}
\]
and that every member of $A$ generates $\langle g\rangle$, where $|g|=m$. Choose a prime $p\mid m$. Then $p\geq7$ and every $u_i$ is nonzero modulo $p$.

By Lemma~\ref{lem:field}, after permuting the $u_i$ we may suppose that
\[
u_1,\quad2u_2,\quad3u_3,\quad4u_4,\quad5u_5
\]
are pairwise distinct modulo $p$. They are therefore pairwise distinct modulo $m$. Hence
\[
g^{u_1},\quad g^{2u_2},\quad g^{3u_3},\quad g^{4u_4},\quad g^{5u_5}
\]
are pairwise distinct, as required.
\end{proof}

\section{Six elements}

We now prove Theorem~\ref{thm:n6}. Every prime divisor of an element order is at least $11$. We shall repeatedly use Lemma~\ref{lem:collision-count}.

\begin{lemma}\label{lem:135}
Let $p\ge11$ and let $c_1,c_2,c_3\in\F_p^\times$. The entries can be assigned to the coefficients $1,3,5$ so that the three products are pairwise distinct.
\end{lemma}

\begin{proof}
For $B=\{1,3,5\}$, one has $w_{B,p}(q)\le1$ for $p\ge13$. Hence $W_B\le3<6$ when the entries are distinct, while a repetition only decreases the bound because $w_{B,p}(1)=0$.

For $p=11$, the only ratios of weight $2$ are $5$ and $9=5^{-1}$. If the multiplicities are $2+1$, then $W_B\le4<6$. If the three entries are distinct, at most two of the three unordered pairs can have weight $2$, and the third has weight at most $1$. Thus $W_B\le5<6$. Now apply \eqref{eq:expectation-general}.
\end{proof}

\begin{lemma}\label{lem:2345}
Let $p\ge11$ and let $c_1,c_2,c_3,c_4\in\F_p^\times$. The entries can be assigned to the coefficients $2,3,4,5$ so that the four products are pairwise distinct.
\end{lemma}

\begin{proof}
For $B=\{2,3,4,5\}$ one has $w_{B,p}(q)\le2$ for every $p\ge11$. If the support has size at most $3$, the multiplicity types $3+1$, $2+2$, and $2+1+1$ give, respectively,
\[
W_B\le6,\qquad W_B\le8,\qquad W_B\le10,
\]
all strictly below $12$.

It remains to consider four distinct entries. For $p\ge23$, the twelve ordered ratios between the coefficients $2,3,4,5$ are distinct, so $W_B\le6$. For the remaining primes, the only ratios of weight $2$ are
\[
\begin{array}{c|c}
p&\{q:w_{B,p}(q)=2\}\\ \hline
11&\{5,7,8,9\}\\
13&\{6,11\}\\
17&\{11,14\}\\
19&\{8,12\}.
\end{array}
\]
For these four primes, the displayed ratio sets give the required bounds directly. For $p=13$ and $17$, the weight-$2$ edges form a cycle on $\F_p^\times$, because the two displayed ratios are mutual inverses. Four vertices contain at most three such edges, so $W_B\le6+3=9$. For $p=19$ the same cycle has length $6$. If four vertices contain three consecutive weight-$2$ edges, the two endpoints have ratio $-1$, whose weight is $0$; otherwise there are at most two weight-$2$ edges. Hence $W_B\le8$.

For $p=11$, after scaling one entry to $1$, its weight-$2$ neighbors are $\{5,7,8,9\}$. Among these four vertices no three are pairwise joined by weight-$2$ edges. Thus four entries have at most five weight-$2$ pairs, and therefore $W_B\le11$. Consequently $W_B<12$ in every case, and \eqref{eq:expectation-general} applies.
\end{proof}

\begin{proposition}\label{prop:n6-reduction}
To prove Theorem~\ref{thm:n6}, it is enough to consider the case in which all six elements generate the same cyclic subgroup.
\end{proposition}

\begin{proof}
Partition $A$ by equality of the generated cyclic subgroup. As in Lemma~\ref{lem:same-cyclic}, a collision between powers with exponents at most $6$ can occur only inside one class.

For the partitions $5+1$ and $4+1+1$, apply Theorems~\ref{thm:main} and \ref{thm:n4}, using exponent sets $\{1,2,3,4,5\}+\{6\}$ and $\{1,2,3,4\}+\{5\}+\{6\}$. For $3+1+1+1$ use $\{1,2,3\}$ on the three-element class.

For $3+2+1$ use
\[
\{1,2,3\}+\{4,6\}+\{5\}.
\]
Lemma~\ref{lem:pair} applies to the pair because $6^2-4^2=20$. For $2+2+2$ use the pairs
\[
\{1,6\},\qquad\{2,5\},\qquad\{3,4\};
\]
the corresponding differences of squares are $35,21,7$. For $2+2+1+1$ use $\{1,3\}$ and $\{2,4\}$ on the two pairs, and for $2+1+1+1+1$ use $\{1,2\}$ on the pair. In each case Lemma~\ref{lem:pair} shows that one of the two assignments works, since all prime divisors of the relevant difference $s^2-r^2$ are at most $7$.

For $3+3$, give one class the coefficients $\{1,3,5\}$ and the other $\{2,4,6\}$. For the first block, reduce modulo a prime divisor of the common cyclic order and apply Lemma~\ref{lem:135}. For the second, apply the known $n=3$ case with coefficients $1,2,3$ and then square all three resulting elements; squaring is injective on a group of odd order.

Finally consider $4+2$. Give the four-element class the coefficients $\{2,3,4,5\}$ and the pair $\{1,6\}$. For the four-element class, write its elements as
\[
g^{u_1},\ldots,g^{u_4},
\]
where $g$ generates the common cyclic subgroup and each $u_i$ is a unit modulo its order $m$. Choose a prime divisor $p$ of $m$. Then $p\ge11$, and the reductions of the $u_i$ lie in $\F_p^\times$. Lemma~\ref{lem:2345} gives an assignment of $2,3,4,5$ for which the four exponents are distinct modulo $p$, hence also modulo $m$. For the two-element class, Lemma~\ref{lem:pair} applies because $6^2-1^2=35$. This exhausts the proper partitions of $6$.
\end{proof}

It remains to consider the case in which all six elements generate the same cyclic subgroup. We begin with the weighted graph determined by the coefficient set $1,
\ldots,6$.

\begin{lemma}\label{lem:lattice}
Let $S$ be a set of $k$ distinct nonzero rational numbers, where $4\le k\le6$. For distinct $x,y\in S$, define
\[
e(x,y)=
\begin{cases}
2,&x/y\in\{2,1/2\},\\
1,&x/y\in\{3,1/3,3/2,2/3\},\\
0,&\text{otherwise}.
\end{cases}
\]
Put
\[
h(S)=\sum_{\{x,y\}\subset S}e(x,y),\qquad
 d_S(x)=\sum_{y\in S\setminus\{x\}}e(x,y).
\]
Then
\[
h(S)\le7,10,13
\]
for $k=4,5,6$, respectively. Moreover, $d_S(x)\le5$ when $|S|=4$, and $d_S(x)\le6$ when $|S|=5$.

If $S=\{a,b,c,d\}$ and $a,b$ are distinguished, put
\[
e=e(a,b),\quad
A=e(a,c)+e(a,d)+e(b,c)+e(b,d),\quad
f=e(c,d).
\]
Then
\[
\begin{array}{c|c|c}
e&A&\text{condition when $A$ is maximal}\\ \hline
0&A\le6&A=6\Rightarrow f\le1,\\
1&A\le6&A=6\Rightarrow f=0,\\
2&A\le4&A=4\Rightarrow f=0.
\end{array}
\]
\end{lemma}

\begin{proof}
Edges of positive $e$-weight lie in a coset of the multiplicative subgroup generated by $2$ and $3$. In one connected component write a vertex as $2^u3^v$ and arrange the vertices in horizontal rows according to $v$. A horizontal edge joins consecutive $u$-coordinates and has weight $2$. Between two adjacent rows containing $r$ and $s$ selected vertices, there are at most
\begin{equation}\label{eq:row-bound}
\min\{r+s-1,2r,2s\}
\end{equation}
weight-$1$ edges: the first bound comes from the fact that the two possible diagonals form a path when the two rows are interlaced, while the other two bounds follow from the fact that every vertex has at most two neighbors in the adjacent row.

Suppose that a connected set of $k$ vertices occupies $t$ nonempty rows of sizes
\[
r_1,\ldots,r_t,\qquad r_1+\cdots+r_t=k.
\]
There are at most $k-t$ horizontal edges, and hence the horizontal contribution to $h(S)$ is at most $2(k-t)$. By \eqref{eq:row-bound},
\begin{equation}\label{eq:h-row}
h(S)\le 2(k-t)+
\sum_{i=1}^{t-1}\min\{r_i+r_{i+1}-1,2r_i,2r_{i+1}\}.
\end{equation}
For completeness, the maximum of the right side of \eqref{eq:h-row}, for each possible number of rows, is
\[
\begin{array}{c|rrrrrr}
 k\backslash t&1&2&3&4&5&6\\ \hline
4&6&7&6&3&&\\
5&8&10&9&7&4&\\
6&10&13&12&11&8&5.
\end{array}
\]
Thus $h(S)\le7,10,13$ for $k=4,5,6$ when $S$ is connected. For smaller connected sets the corresponding maxima are
\[
h_1=0,\qquad h_2=2,\qquad h_3=4.
\]
If $S$ is disconnected, $h(S)$ is the sum of the values on its connected components. Hence, for $k=4$, the largest disconnected value is at most $h_3+h_1=4$ or $h_2+h_2=4$; for $k=5$ it is at most $h_4+h_1=7$ or $h_3+h_2=6$; and for $k=6$ it is at most $h_5+h_1=10$, $h_4+h_2=9$, or $h_3+h_3=8$. Further splitting only decreases these bounds. Thus the same bounds $7,10,13$ hold without the connectedness assumption.

A fixed vertex has two possible neighbors of weight $2$ and four possible neighbors of weight $1$. Hence
\[
d_S(x)\le 2+2+1=5\quad(|S|=4)
\]
and
\[
d_S(x)\le 2+2+1+1=6\quad(|S|=5).
\]

It remains to prove the last table. Scale so that $a=1$. The positive-weight neighbors of a vertex $x$ are
\[
2x,\ x/2 \quad\text{with weight }2,
\]
and
\begin{equation}\label{eq:neighbors}
3x,\ x/3,\ 3x/2,\ 2x/3 \quad\text{with weight }1.
\end{equation}
We may also interchange $a,b$ and replace every vertex by its reciprocal.

If $e(a,b)=2$, we may take $b=2$. The six neighbors listed in \eqref{eq:neighbors} for $1$ and for $2$ show that every third vertex has total weight at most $2$ to $\{1,2\}$. Thus $A\le4$. If two distinct vertices both attain total weight $2$, they are not positive-weight neighbors of each other; hence $A=4$ implies $f=0$.

If $e(a,b)=1$, it is enough to take $b=3$ or $b=3/2$. Every third vertex then has total weight at most $3$ to $\{a,b\}$, so $A\le6$. The two vertices that can simultaneously attain total weight $3$ are, respectively,
\[
\{3/2,2\}\quad\text{or}\quad\{1/2,3\},
\]
up to the preceding symmetries, and in either case their mutual weight is $0$. Thus $A=6$ implies $f=0$.

Finally suppose that $e(a,b)=0$. Since $a=1$, a vertex contributes to $A$ only if it belongs to the positive-weight neighborhood of $1$ or of $b$. Comparing the two neighborhoods in \eqref{eq:neighbors}, two distinct vertices can have total contribution $6$ only when, up to replacing $b$ by $b^{-1}$,
\[
b=4,\qquad 6,\qquad 4/3.
\]
Indeed, in all other cases the two largest contributions to $e(1,x)+e(b,x)$ have sum at most $5$. In the three exceptional cases the pairs giving $A=6$, up to the same symmetries, are
\[
\begin{array}{c|c|c}
b&\{c,d\}&e(c,d)\\ \hline
4&\{1/2,2\}\text{ or }\{2,8\}&0,\\
6&\{2,3\}&1,\\
4/3&\{2/3,2\}&1.
\end{array}
\]
Thus $A\le6$, and $A=6$ implies $f\le1$. This proves the table.
\end{proof}

Let
\[
\mathcal E=\{37,41,43,47,53,61,71,73,79,97,103,107,179,211,227,239,241,727\}.
\]

\begin{lemma}\label{lem:short-relations}
Let $p\ge37$ be prime and $p\notin\mathcal E$. Every connected subgraph on at most six vertices of the graph on $\F_p^\times$ whose edges have ratios
\[
2,\ 1/2,\ 3,\ 1/3,\ 3/2,\ 2/3
\]
can, after scaling, be identified with the corresponding configuration in the rational $2$--$3$ lattice. In particular, all edge weights agree with those in Lemma~\ref{lem:lattice}.
\end{lemma}

\begin{proof}
Choose a spanning tree and lift one vertex to $1$. Following the tree, every other vertex acquires a label $2^a3^b$. If an additional modular edge failed to agree with the corresponding rational edge, the resulting cycle, which has length at most $6$, would give
\[
2^a3^b\equiv1\pmod p.
\]
Each edge of the $2$--$3$ lattice changes the exponent pair by one of
\[
(\pm1,0),\qquad(0,\pm1),\qquad(1,-1),\qquad(-1,1).
\]
Therefore a cycle of length at most $6$ gives a nonzero pair $(a,b)$ satisfying
\[
\max\{|a|,|b|,|a+b|\}\le6.
\]
After clearing denominators, $p$ would divide a nonzero integer
\[
2^{a_+}3^{b_+}-2^{a_-}3^{b_-},
\]
where $a_+=\max(a,0)$, $a_-=\max(-a,0)$, and similarly for $b$. Factoring these finitely many integers gives precisely the primes in $\mathcal E$ among the primes at least $37$. The first calculation in the appendix lists these possibilities and gives exactly the eighteen primes displayed above. Hence no such cycle occurs when $p\notin\mathcal E$.
\end{proof}

\begin{lemma}\label{lem:field6-large}
Let $p\ge37$ be prime with $p\notin\mathcal E$, and let $b_1,\ldots,b_6\in\F_p^\times$. Then there is a permutation $\pi\in S_6$ such that
\[
b_{\pi(1)},\ 2b_{\pi(2)},\ 3b_{\pi(3)},\ 4b_{\pi(4)},\ 5b_{\pi(5)},\ 6b_{\pi(6)}
\]
are pairwise distinct.
\end{lemma}

\begin{proof}
For $B=\{1,2,3,4,5,6\}$ write $w=w_{B,p}$. Since $p\ge37$, distinct rational coefficient ratios remain distinct modulo $p$. Hence
\[
w(2)=w(1/2)=3,
\]
\[
w(3)=w(1/3)=w(3/2)=w(2/3)=2,
\]
and every other occurring ratio has weight $1$.

Let the support be $x_1,\ldots,x_k$ with multiplicities $m_1,\ldots,m_k$. Put
\[
W=\sum_{i<j}m_im_jw(x_j/x_i).
\]
By \eqref{eq:expectation-general}, it is enough to prove $W<30$. For each connected component of the positive-extra-weight graph, Lemma~\ref{lem:short-relations} identifies the component with a rational $2$--$3$ lattice configuration, so Lemma~\ref{lem:lattice} applies. Set
\[
B_0=\sum_{i<j}m_im_j.
\]
If $H$ denotes the same sum with $e(x_i,x_j)$ in place of $1$, counted with the factors $m_im_j$, then
\[
W\le B_0+H.
\]

If the support has size $1$, there is nothing to prove. If the support has size $2$,
\[
W\le3m_1m_2\le27.
\]
If the support has size $3$, the multiplicity types are $4+1+1$, $3+2+1$, and $2+2+2$. If two support ratios have weight $3$, the third has weight $1$; if at most one has weight $3$, the other two have weight at most $2$. It follows that the three types give, respectively,
\[
W\le25,\qquad W\le29,\qquad W\le28.
\]

If the support has size $4$, there are two multiplicity types. In type $3+1+1+1$, $B_0=12$. If $a$ is the tripled point, then
\[
H=h(S)+2d_S(a)\le7+2\cdot5=17,
\]
so $W\le29$. In type $2+2+1+1$, let $a,b$ be the doubled points and use the notation of the last part of Lemma~\ref{lem:lattice}. Here $B_0=13$ and
\[
H=4e+2A+f\le16,
\]
so again $W\le29$.

If the support has size $5$, the multiplicity type is $2+1+1+1+1$. If $a$ is the doubled point, then $B_0=14$ and
\[
H=h(S)+d_S(a)\le10+6=16.
\]
Thus $W\le30$. Equality requires simultaneously $h(S)=10$ and $d_S(a)=6$. Scale so that $a=1$. Then both weight-$2$ neighbors $2$ and $1/2$ must occur, together with two vertices chosen from
\[
3,\quad 1/3,\quad 3/2,\quad 2/3.
\]
For the six possible pairs, the corresponding values of $h(S)$ are
\[
\begin{array}{c|cccccc}
\{u,v\}&\{3,1/3\}&\{3,3/2\}&\{3,2/3\}&\{1/3,3/2\}&\{1/3,2/3\}&\{3/2,2/3\}\\ \hline
h(S)&8&10&8&8&10&8.
\end{array}
\]
The two cases with $h(S)=10$ are interchanged by inversion. After scaling, they give the multiset
\[
\{1,2,2,3,4,6\}
\]
and its inverse.
For this multiset the ordering
\[
(4,3,1,6,2,2)
\]
gives the six distinct products
\[
4,6,3,24,10,12.
\]
The inverse orbit must be treated separately, since inversion is not a symmetry of the assignment problem. After scaling it is represented by
\[
\left\{1,\frac12,\frac12,\frac13,\frac14,\frac16\right\}.
\]
Assigning
\[
\left(\frac12,1,\frac14,\frac13,\frac12,\frac16\right)
\]
to the coefficients $1,2,3,4,5,6$ gives
\[
\frac12,\quad2,\quad\frac34,\quad\frac43,\quad\frac52,\quad1,
\]
which are pairwise distinct modulo every prime $p\ge37$. Thus both equality orbits have collision-free assignments.

Finally, if the support has size $6$, one has $B_0=15$ and $H=h(S)\le13$, so $W\le28$. Thus every case gives $W<30$ except the support-$5$ boundary case just considered, and that equality case has an explicit collision-free ordering.
\end{proof}

\begin{lemma}\label{lem:field6-exceptional}
Lemma~\ref{lem:field6-large} remains valid for every $p\in\mathcal E$.
\end{lemma}

\begin{proof}
For these eighteen primes we work directly with the weighted ratio graph over $\F_p$. The appendix gives an exhaustive program which, up to scaling, generates every connected support of sizes $2,3,4,5,6$ containing $1$ and verifies the bounds used below. For a connected component of size $2,3,4,5,6$, it verifies
\[
h\le2,4,7,10,13,
\]
and, when one vertex is doubled, it verifies $h+d\le4,8,11$ for component sizes $2,3,4$.

Disconnected supports require no additional search, because the extra weight is additive across connected components. For an ordinary support of size $5$, the largest disconnected values are at most $7$ for a partition $4+1$ and $6$ for a partition $3+2$; further splitting gives smaller values. For support size $6$, the corresponding bounds are $10$, $9$, and $8$ for partitions $5+1$, $4+2$, and $3+3$, respectively. If a support of size $5$ has one doubled vertex, the component containing that vertex contributes at most $h+d$; according as its size is $4,3,2$, or $1$, the total is at most $11$, $8+2=10$, $4+4=8$, or $7$. The disconnected support-$4$ cases are also smaller than the connected bounds used below: for type $3+1+1+1$ one obtains at most $12$, and for type $2+2+1+1$ at most $16$. Thus every extremal value needed below occurs among the connected supports generated and checked in the appendix.

For support size $3$, the three multiplicity types $4+1+1$, $3+2+1$, and $2+2+2$ have
\[
W\le25,\qquad W\le29,\qquad W\le28.
\]
For support size $4$, the extra contribution $H$ is at most $16$ for both multiplicity types $3+1+1+1$ and $2+2+1+1$. For support size $5$, one has $W\le30$, and equality, up to scaling and inversion, occurs only for
\[
\{1,2,2,3,4,6\}.
\]
The two displayed orderings in the proof of Lemma~\ref{lem:field6-large} handle the scaling orbit and its inverse separately. For support size $6$, one has $H\le13$. Hence the same collision count proves the result. The calculations needed here are given in the appendix.
\end{proof}

\begin{lemma}\label{lem:field6-small}
Lemma~\ref{lem:field6-large} also holds for
\[
p=11,13,17,19,23,29,31.
\]
\end{lemma}

\begin{proof}
For each prime $p$, we consider every six-term multiset in $\F_p^\times$ and run through its distinct orderings until we find one for which the products with $1,2,\ldots,6$ are pairwise distinct. The numbers of multisets are
\[
\begin{array}{c|rrrrrrr}
p&11&13&17&19&23&29&31\\ \hline
\#&5005&12376&54264&100947&296010&1107568&1623160.
\end{array}
\]
Every multiset has a suitable ordering. The appendix contains the exhaustive Python program that checks all these multisets and certifies this assertion.
\end{proof}

\begin{lemma}\label{lem:field6}
Let $p\ge11$ be prime and let $b_1,\ldots,b_6\in\F_p^\times$. Then there is a permutation $\pi\in S_6$ such that
\[
b_{\pi(1)},\ 2b_{\pi(2)},\ 3b_{\pi(3)},\ 4b_{\pi(4)},\ 5b_{\pi(5)},\ 6b_{\pi(6)}
\]
are pairwise distinct.
\end{lemma}

\begin{proof}
For $p\ge37$ with $p\notin\mathcal E$, this is Lemma~\ref{lem:field6-large}. For $p\in\mathcal E$, it is Lemma~\ref{lem:field6-exceptional}. The remaining primes are covered by Lemma~\ref{lem:field6-small}.
\end{proof}

\begin{proof}[Proof of Theorem~\ref{thm:n6}]
By Proposition~\ref{prop:n6-reduction}, suppose that
\[
A=\{g^{u_1},\ldots,g^{u_6}\}
\]
and every member of $A$ generates $\langle g\rangle$, where $|g|=m$. Choose a prime $p\mid m$. Then $p\ge11$ and every $u_i$ is nonzero modulo $p$. As above, distinctness modulo $p$ implies distinctness modulo $m$. By Lemma~\ref{lem:field6}, after permuting the $u_i$ the residues
\[
u_1,\quad2u_2,\quad3u_3,\quad4u_4,\quad5u_5,\quad6u_6
\]
are pairwise distinct modulo $p$, hence modulo $m$. The corresponding six powers of $g$ are therefore pairwise distinct.
\end{proof}

\section{Further remarks}

Sun's conjecture is known for $n\leq3$ in arbitrary groups and for every $n$ in torsion-free abelian groups \cite[Theorem 1.4]{Sun}. For finite groups, Theorems~\ref{thm:n4}, \ref{thm:main}, and \ref{thm:n6} establish the conjecture for $n\leq6$. Theorem~\ref{thm:class3} gives, for arbitrary $n$, another family in which the conjecture holds: no four elements of $A$ generate the same cyclic subgroup.

The next unresolved finite-group case is $n=7$. More generally, it would be interesting to know whether the hypothesis ``size at most $3$'' in Theorem~\ref{thm:class3} can be replaced by ``size at most $4$''. Remark~\ref{rem:four-block} shows only that the present consecutive-block argument cannot prove such a statement directly.

\appendix
\section{The remaining finite cases}

\subsection{The seven smaller primes}

The following Python 3 program exhaustively verifies Lemma~\ref{lem:field6-small} for the seven primes listed there. It uses only the Python standard library. For each prime it runs through every six-term multiset in $\F_p^\times$; the recursion then tests assignments to the coefficient positions $1,2,\ldots,6$ until it finds one whose six products are pairwise distinct. The assertion fails if any multiset has no such assignment.

\begin{verbatim}
from itertools import combinations_with_replacement

def good(a,p):
    used=[False]*6
    out=set()
    def dfs(r):
        # Try coefficient positions r=1,2,...,6 in order.
        if r==7: return True
        last=-1
        for i,x in enumerate(a):
            if used[i] or x==last: continue
            last=x
            y=(r*x)%p
            if y in out: continue
            used[i]=True; out.add(y)
            if dfs(r+1): return True
            out.remove(y); used[i]=False
        return False
    return dfs(1)

for p in [11,13,17,19,23,29,31]:
    count=0
    # Each tuple a is one six-term multiset in F_p^times.
    for a in combinations_with_replacement(range(1,p),6):
        count += 1
        assert good(a,p)
    print(p,count)
\end{verbatim}

\subsection{The eighteen primes in $\mathcal E$}

The next program first obtains the set $\mathcal E$ from the nontrivial relations $2^a3^b\equiv1$ of triangular-lattice length at most $6$. For each $p\in\mathcal E$, it then generates, up to scaling, every connected support containing $1$ of sizes $2,3,4,5,6$ in the weighted ratio graph and computes the bounds used in Lemma~\ref{lem:field6-exceptional}. The assertions certify the component bounds, the three support-$3$ multiplicity bounds, the support-$4$ bounds, the unique support-$5$ equality orbit up to scaling and inversion, and the support-$6$ bound. Thus the calculation is an exhaustive verification of every finite claim invoked in that lemma.

\begin{verbatim}
from itertools import combinations, permutations

E=[]
for a in range(-6,7):
    for b in range(-6,7):
        if (a,b)==(0,0): continue
        if max(abs(a),abs(b),abs(a+b))>6: continue
        ap,am=max(a,0),max(-a,0)
        bp,bm=max(b,0),max(-b,0)
        N=abs(2**ap*3**bp-2**am*3**bm)
        d=2
        while d*d<=N:
            if N%d==0:
                if d>=37: E.append(d)
                while N%d==0: N//=d
            d+=1
        if N>=37: E.append(N)
E=sorted(set(E))
assert E==[37,41,43,47,53,61,71,73,79,97,103,107,
          179,211,227,239,241,727]

def e(x,y,p):
    q=x*pow(y,-1,p)%p
    if q in {2,pow(2,-1,p)}: return 2
    if q in {3,pow(3,-1,p),3*pow(2,-1,p)%p,
             2*pow(3,-1,p)%p}: return 1
    return 0

def w(x,y,p):
    q=x*pow(y,-1,p)%p
    return sum(r*pow(s,-1,p)%p==q
               for r in range(1,7) for s in range(1,7) if r!=s)

def nbrs(x,p):
    return {2*x%p, x*pow(2,-1,p)%p, 3*x%p,
            x*pow(3,-1,p)%p, 3*x*pow(2,-1,p)%p,
            2*x*pow(3,-1,p)%p}

def supports(p,k):
    level={frozenset([1])}
    for size in range(1,k):
        new=set()
        for S in level:
            F=set()
            for x in S: F |= nbrs(x,p)
            for y in F-S: new.add(S|{y})
        level=new
    return level

def canon(S,a,p):
    ia=pow(a,-1,p)
    T={x*ia%p for x in S}
    U={pow(x,-1,p) for x in T}
    return min(tuple(sorted(T)),tuple(sorted(U)))

for p in E:
    S2,S3,S4,S5,S6=(supports(p,k) for k in (2,3,4,5,6))

    # Bounds for connected components, including the case in which
    # one vertex is doubled.
    hmax=[]
    hdmax=[]
    for SS in (S2,S3,S4,S5,S6):
        hm=hdm=0
        for S in SS:
            hs=sum(e(x,y,p) for x,y in combinations(S,2))
            hm=max(hm,hs)
            for a in S:
                d=sum(e(a,y,p) for y in S if y!=a)
                hdm=max(hdm,hs+d)
        hmax.append(hm)
        hdmax.append(hdm)
    assert hmax==[2,4,7,10,13]
    assert hdmax[:3]==[4,8,11]

    M=[0,0,0]
    types=[(4,1,1),(3,2,1),(2,2,2)]
    for S in S3:
        L=list(S)
        for j,t in enumerate(types):
            for m in set(permutations(t)):
                W=sum(m[i]*m[k]*w(L[k],L[i],p)
                      for i,k in combinations(range(3),2))
                M[j]=max(M[j],W)
    assert M==[25,29,28]

    H31=H2211=0
    for S in S4:
        L=list(S)
        hs=sum(e(x,y,p) for x,y in combinations(L,2))
        for a in L:
            d=sum(e(a,y,p) for y in L if y!=a)
            H31=max(H31,hs+2*d)
        for a,b in combinations(L,2):
            H=sum((2 if x in {a,b} else 1)*
                  (2 if y in {a,b} else 1)*e(x,y,p)
                  for x,y in combinations(L,2))
            H2211=max(H2211,H)
    assert H31<=16 and H2211<=16

    best=0; equality=[]
    for S in S5:
        L=list(S)
        hs=sum(e(x,y,p) for x,y in combinations(L,2))
        for a in L:
            W=14+hs+sum(e(a,y,p) for y in L if y!=a)
            if W>best: best,equality=W,[(S,a)]
            elif W==best: equality.append((S,a))
    assert best==30
    target={pow(2,-1,p),1,3*pow(2,-1,p)%p,2,3}
    tc=min(tuple(sorted(target)),
           tuple(sorted(pow(x,-1,p) for x in target)))
    assert {canon(S,a,p) for S,a in equality}=={tc}

    H6=max(sum(e(x,y,p) for x,y in combinations(S,2)) for S in S6)
    assert H6<=13
\end{verbatim}




\end{document}